\documentclass[12pt,high]{amsart}
\usepackage{amsmath,amsthm,amsfonts,amssymb,mathrsfs}
\date{\today}

\usepackage[cp1251]{inputenc}         

\usepackage[ukrainian]{babel} 
\usepackage{amsmath,amsthm,amsfonts,amssymb,mathrsfs}
\usepackage{eucal}

\usepackage{color}

\usepackage{amssymb,cite}
\usepackage[colorlinks,plainpages,citecolor=magenta, linkcolor=blue, backref]{hyperref}

\usepackage{hyperref}

\newtheorem{theorem}{Теорема}

\newtheorem{proposition}{Твердження}
\newtheorem{corollary}{Наслiдок}
\newtheorem{lemma}{Лема}
\theoremstyle{definition}
\newtheorem{remark}{Зауваження}

\begin{document}

\title[Про ендоморфiзми бiциклiчно\"{\i} напiвгрупи та розширено\"{\i} бiциклiчно\"{\i} напiвгрупи]{Про ендоморфiзми бiциклiчно\"{\i} напiвгрупи та розширено\"{\i} бiциклiчно\"{\i} напiвгрупи}

\author[Олег~Гутік, Оксана Прохоренкова, Діана Сех]{Олег~Гутік, Оксана Прохоренкова, Діана Сех}
\address{Механіко-математичний факультет, Львівський національний університет ім. Івана Франка, Університецька 1, Львів, 79000, Україна}
\email{oleg.gutik@lnu.edu.ua, 
oksana.prokhorenkova@lnu.edu.ua,  diana.sekh@lnu.edu.ua}



\begin{abstract}
Доведено, що напівгрупи $\mathrm{\mathbf{End}}(\boldsymbol{B}_{\omega})$ та $\mathrm{\mathbf{End}}(\boldsymbol{B}_{\mathbb{Z}})$ ендоморфізмів  бі\-цик\-ліч\-ної напівгрупи $\boldsymbol{B}_{\omega}$ та ендоморфізмів розширеної біциклічної напівгрупи $\boldsymbol{B}_{\mathbb{Z}}$ ізоморфні напівпрямим добуткам $(\omega,+)\rtimes_\varphi(\omega,*)$ і $\mathbb{Z}(+)\rtimes_\varphi(\omega,*)$, відповідно.

\bigskip
\noindent
\emph{Oleg Gutik, Oksana Prokhorenkova, Diana Sekh, \textbf{On endomorphisms of the bicyclic semigroup and the extended bicyclic semigroup}.}

\smallskip
\noindent
It is proved that the semigroups $\mathrm{\mathbf{End}}(\boldsymbol{B}_{\omega})$ and $\mathrm{\mathbf{End}}(\boldsymbol{B}_{\mathbb{Z}})$ of the endomorphisms of the bicyclic semigroup $\boldsymbol{B}_{\omega}$ and the endomorphisms of the extended bicyclic semigroup $\boldsymbol{B}_{\mathbb{Z}}$ are isomorphic to the semidirect products $(\omega,+)\rtimes_\varphi(\omega,*)$ and $\mathbb{Z}(+)\rtimes_\varphi(\omega,*)$, respectively.

\end{abstract}

\maketitle


\section{\textbf{Вступ}}\label{section-1}

Ми користуватимемося термінологією з \cite{Clifford-Preston-1961, Clifford-Preston-1967, Lawson-1998, Petrich1984}.
Надалі у тексті множину невід'ємних цілих чисел  позначатимемо через $\omega$, множину цілих чисел через $\mathbb{Z}$, і адитивну групу цілих чисел через $\mathbb{Z}(+)$. Надалі, якщо $f\colon X\to Y$~--- відображення, то через $(x)f$ будемо позначати образ елемента $x\in X$ стосовно $f$.

\smallskip

Якщо $S$~--- напівгрупа, то її підмножина ідемпотентів позначається через $E(S)$.  На\-пів\-гру\-па $S$ називається \emph{інверсною}, якщо для довільного її елемента $x$ існує єдиний елемент $x^{-1}\in S$ такий, що $xx^{-1}x=x$ та $x^{-1}xx^{-1}=x^{-1}$ \cite{Vagner-1952, Petrich1984}. В інверсній напівгрупі $S$ вище означений елемент $x^{-1}$ називається \emph{інверсним до} $x$. \emph{В'язка}~--- це напівгрупа ідемпотентів, а \emph{напівґратка}~--- це комутативна в'язка.

\smallskip

Через $(\omega,+)$ і $(\omega,*)$ позначатимемо адитивну та мультиплікативну напівгрупи невід'ємних цілих чисел, відповідно.

\smallskip

Якщо $S$~--- напівгрупа, то на $E(S)$ визначено частковий порядок:
$
e\preccurlyeq f
$   тоді і лише тоді, коли
$ef=fe=e$.
Так означений частковий порядок на $E(S)$ називається \emph{при\-род\-ним}.

\smallskip

Означимо відношення $\preccurlyeq$ на інверсній напівгрупі $S$ так:
$
    s\preccurlyeq t
$
тоді і лише тоді, коли $s=te$.
для деякого ідемпотента $e\in S$. Так означений частковий порядок назива\-єть\-ся \emph{при\-род\-ним част\-ковим порядком} на інверсній напівгрупі $S$~\cite{Vagner-1952}. Оче\-вид\-но, що звуження природного часткового порядку $\preccurlyeq$ на інверсній напівгрупі $S$ на її в'язку $E(S)$ є при\-род\-ним частковим порядком на $E(S)$.

\smallskip

Нехай $S$ i $T$~--- напівгрупи. Відображення $\mathfrak{h}\colon S\to T$ називається \emph{гомомор\-фіз\-мом}, якщо $(x\cdot y)\mathfrak{h}=(x)\mathfrak{h}\cdot (y)\mathfrak{h}$ для довільних $x,y\in S$ \cite{Clifford-Preston-1961}. Якщо ж $S$ i $T$~--- моноїди з одиницями $1_S$ i $1_T$, відповідно, то гомоморфізм $\mathfrak{h}\colon S\to T$ такий, що $(1_S)\mathfrak{h}=1_T$, будемо називати \emph{гомоморфізмом моноїдів}. Також гомоморфізм $\mathfrak{h}\colon S\to S$ на\-зи\-ва\-єть\-ся \emph{ендоморфізмом} напівгрупи $S$,  ізоморфізм $\mathfrak{i}\colon S\to S$ на\-зи\-ва\-єть\-ся \emph{автомор\-фіз\-мом} напівгрупи $S$, а якщо $S$~--- моноїд з одиницею $1_S$ і $\mathfrak{h}\colon S\to S$ --- гомоморфізм такий, що $(1_S)\mathfrak{h}=1_S$, то будемо говорити, що $\mathfrak{h}$ --- \emph{ендоморфізм} моноїда $S$. Очевидно, що композиція двох ендоморфізмів напівгрупи (моноїда) $S$ є ендоморфізмом напівгрупи (моноїда) $S$, а композиція двох автоморфізмів напівгрупи (моноїда) $S$ є автоморфізмом напівгрупи (моноїда) $S$. Отож множина усіх ендомор\-фіз\-мів $\mathrm{\mathbf{End}}(S)$ напівгрупи $S$ стосовно операції композиції відображень є моноїдом, а множина усіх авто\-мор\-фіз\-мів $\mathrm{\mathbf{Aut}}(S)$ напівгрупи $S$ стосовно операції композиції відображень є групою, і очевидно є групою одиниць напівгрупи $\mathrm{\mathbf{End}}(S)$. Якщо $S=S^1$~--- моноїд, то надалі напівгрупу ендоморфізмів моноїда $S$ (як моноїда) будемо позначати через $\mathrm{\mathbf{End}}^1(S)$, а напівгрупу ендоморфізмів напівгрупи $S$ (як напівгрупи) будемо позначати через $\mathrm{\mathbf{End}}(S)$.

\smallskip

Нагадаємо (див. \cite[\S1.12]{Clifford-Preston-1961}), що \emph{біциклічною напівгрупою} (або \emph{біциклічним моноїдом}) ${\mathscr{C}}(p,q)$ називається напівгрупа з одиницею, породжена двоелементною мно\-жи\-ною $\{p,q\}$ і визначена одним  співвідношенням $pq=1$. Біциклічна на\-пів\-група відіграє важливу роль у теорії на\-півгруп. Зокрема, класична теорема О.~Ан\-дерсена \cite{Andersen-1952}  стверд\-жує, що {($0$-)}прос\-та напівгрупа з (ненульовим) ідем\-по\-тен\-том є цілком {($0$-)}прос\-тою тоді і лише тоді, коли вона не містить ізоморфну копію бі\-циклічної напівгрупи.

\begin{remark}\label{remark-10}
Легко бачити, що біциклічний моноїд ${\mathscr{C}}(p,q)$ ізоморфний напівгрупі, заданій на множині $\boldsymbol{B}_{\omega}=\omega\times\omega$ з напівгруповою операцією
\begin{equation}\label{eq-1.1}
\begin{split}
    (i_1,j_1)\cdot(i_2,j_2)&=(i_1+i_2-\min\{j_1,i_2\},j_1+j_2-\min\{j_1,i_2\})=\\
  &=
\left\{
  \begin{array}{ll}
    (i_1-j_1+i_2,j_2), & \hbox{якщо~} j_1\leqslant i_2;\\
    (i_1,j_1-i_2+j_2), & \hbox{якщо~} j_1\geqslant i_2,
  \end{array}
\right.
\end{split}
\end{equation}
причому цей ізоморфізм визначається за формулою $(q^ip^j)\mathfrak{i}=(i,j)$, $i,j\in\omega$.
\end{remark}

Піднапівгрупи біциклічного моноїда $\boldsymbol{B}_{\omega}$ описано в багатьох працях (див. \cite{Hovsepyan-2017, Hovsepyan-2020, Descalco-Ruskuc-2005, Jin-Guo-You-2017, Makanjuola-Umar-1997, Ovsepyan-2014}). Узагальнення відношень Ґріна, сумісні часткові порядки, авто\-мор\-фіз\-ми, напівавтоморфізми та конгруенції на $\boldsymbol{B}_{\omega}$ вивчали в \cite{Descalco-Higgins-2010, Duchamp-1986, Makanjuola-Umar-1997, McAlister-1999, Nico-1978, Taguchi-Saito-1967}. Інші властивості біциклічного моноїда описано в \cite{Bouvier-Faisant-1970, Clifford-Preston-1961, Lawson-1998, Petrich1984}.

\smallskip

Добре відомо, що кожен неін'єктивний гомоморфний образ біциклічної напівгрупи є циклічною групою (див. \cite[наслідок~1.32]{Clifford-Preston-1961}), а також лише тотожне відображення біциклічної напівгрупи є її автоморфізмом.

\smallskip

Нагадаємо (див. \cite{Warne-1968}), що \emph{розширеною біциклічною напівгрупою} називається множина $\boldsymbol{B}_{\mathbb{Z}}=\mathbb{Z}\times\mathbb{Z}$ з напівгруповою операцією \eqref{eq-1.1}. Відношення Ґріна на розширеній біциклічній напівгрупі описані в \cite{Fihel-Gutik-2011}. Також кожний неін'єктивний гомо\-морф\-ний образ розширеної біциклічної напівгрупи $\boldsymbol{B}_{\mathbb{Z}}$ є циклічною групою \cite{Fihel-Gutik-2011}. Група $\mathrm{\mathbf{Aut}}(\boldsymbol{B}_{\mathbb{Z}})$ автоморфізмів напівгрупи $\boldsymbol{B}_{\mathbb{Z}}$ ізоморфна адитивній групі цілих чисел $\mathbb{Z}(+)$ \cite{Gutik-Maksymyk-2017}.

\smallskip

Варіанти біциклічного моноїда та розширеної біциклічної напівгрупи вивчали в \cite{Givens-Rosin-Linton-2017, Gutik-Maksymyk-2017}.

\smallskip

Ми описуємо ендоморфізми біциклічного моноїда $\boldsymbol{B}_{\omega}$ як моноїда та як напівгрупи. Доведено, що напівгрупа $\mathrm{\mathbf{End}}(\boldsymbol{B}_{\omega})$ ендоморфізмів біциклічної напівгрупи $\boldsymbol{B}_{\omega}$ ізоморфна напівпрямому добутку $(\omega,+)\rtimes_\varphi(\omega,*)$, а також, що напівгрупа $\mathrm{\mathbf{End}}(\boldsymbol{B}_{\mathbb{Z}})$ ендоморфізмів розширеної біциклічної напівгрупи $\boldsymbol{B}_{\mathbb{Z}}$ ізоморфна напівпрямому добутку $\mathbb{Z}(+)\rtimes_\varphi(\omega,*)$.

\section{\textbf{Напівгрупа ендоморфізмів біциклічної напівгрупи}}\label{section-2}

Для довільних $a,k\in\omega$ означимо відображення $\varepsilon_{k[a]}\colon \boldsymbol{B}_{\omega}\to \boldsymbol{B}_{\omega}$ за формулою
\begin{equation}\label{eq-2.1}
  (m,n)\varepsilon_{k[a]}=(km+a,kn+a), \qquad m,n\in \omega.
\end{equation}

Якщо $k=0$, то $(m,n)\varepsilon_{0[a]}=(a,a)$ для всіх $(m,n)\in \boldsymbol{B}_{\omega}$, і відображення $\varepsilon_{0[a]}\colon \boldsymbol{B}_{\omega}\to \boldsymbol{B}_{\omega}$  анулюючий ендоморфізм напівгрупи $\boldsymbol{B}_{\omega}$. Тому надалі вважатимемо, що $k>0$. Тоді для довільних $(i,j),(m,n)\in \boldsymbol{B}_{\omega}$ маємо, що
\begin{align*}
  ((i,j)\cdot(m,n))\varepsilon_{k[a]} &=
  \left\{
  \begin{array}{ll}
    (i-j+m,n)\varepsilon_{k[a]}, & \hbox{якщо~} j<m;\\
    (i,n)\varepsilon_{k[a]},     & \hbox{якщо~} j=m;\\
    (i,j-m+n)\varepsilon_{k[a]}, & \hbox{якщо~} j>m
  \end{array}
\right.=
  \\
   &=
     \left\{
  \begin{array}{ll}
    (k(i-j+m)+a,kn+a), & \hbox{якщо~} j<m;\\
    (ki+a,kn+a),       & \hbox{якщо~} j=m;\\
    (ki+a,k(j-m+n)+a), & \hbox{якщо~} j>m
  \end{array}
\right.
\end{align*}
i оскільки $k>0$, то
\begin{align*}
  (i,j)\varepsilon_{k[a]}&\cdot (m,n)\varepsilon_{k[a]}=(ki+a,kj+a)\cdot (km+a,kn+a)= \\
  &=
  \left\{
  \begin{array}{ll}
    (ki+a-(kj+a)+km+a,kn+a), & \hbox{якщо~} kj+a<km+a;\\
    (ki+a,kn+a),             & \hbox{якщо~} kj+a=km+a;\\
    (ki+a,kj+a-(km+a)+kn+a), & \hbox{якщо~} kj+a>km+a
  \end{array}
\right.= \\
  &=
  \left\{
  \begin{array}{ll}
    (k(i-j+m)+a,kn+a), & \hbox{якщо~} kj<km;\\
    (ki+a,kn+a),             & \hbox{якщо~} kj=km;\\
    (ki+a,k(j-m+n)+a), & \hbox{якщо~} kj>km
  \end{array}
\right.=
  \\
   &=
     \left\{
  \begin{array}{ll}
    (k(i-j+m)+a,kn+a), & \hbox{якщо~} j<m;\\
    (ki+a,kn+a),       & \hbox{якщо~} j=m;\\
    (ki+a,k(j-m+n)+a), & \hbox{якщо~} j>m.
  \end{array}
\right.
\end{align*}

Отже. ми довели таку лему.

\begin{lemma}\label{lemma-2.1}
Для довільних $a,k\in\omega$ відображення $\varepsilon_{k[a]}\colon \boldsymbol{B}_{\omega}\to \boldsymbol{B}_{\omega}$, означене за формулою \eqref{eq-2.1} є ендоморфізмом біциклічної напівгрупи.
\end{lemma}

З леми~\ref{lemma-2.1} випливає наслідок.

\begin{corollary}\label{corollary-2.2}
Для довільного числа $k\in\omega$ відображення $\varepsilon_{k[0]}\colon \boldsymbol{B}_{\omega}\to \boldsymbol{B}_{\omega}$, означене за формулою \eqref{eq-2.1}, де $a=0$, є ендоморфізмом біциклічного моноїда.
\end{corollary}

\begin{lemma}\label{lemma-2.3}
Для довільного ендоморфізму $\varphi\colon \boldsymbol{B}_{\omega}\to \boldsymbol{B}_{\omega}$ біциклічного моноїда існує таке число $k\in\omega$, що $\varphi=\varepsilon_{k[0]}$.
\end{lemma}

\begin{proof}
Оскільки $\varphi\colon\boldsymbol{B}_{\omega}\to\boldsymbol{B}_{\omega}$ --- ендоморфізм біциклічного моноїда, то $(0,0)\varphi=(0,0)$. Тоді для породжуючих елементів $(0,1)$ i $(1,0)$ біциклічного моноїда $\boldsymbol{B}_{\omega}$ існують елементи $(i,j),(m,n)\in \boldsymbol{B}_{\omega}$ такі, що $(0,1)\varphi=(i,j)$ i $(1,0)\varphi=(m,n)$ для деяких $i,j,m,n\in\omega$. Оскільки $\boldsymbol{B}_{\omega}$ --- інверсна напівгрупа та $(0,1)$ i $(1,0)$ --- інверсні елементи в $\boldsymbol{B}_{\omega}$, то
\begin{align*}
  (m,n)& =(1,0)\varphi=\\
    & =((0,1)^{-1})\varphi=\\
    & =((0,1)\varphi)^{-1}=\\
    & =(i,j)^{-1}=(j,i).
\end{align*}
Також з рівності $(0,1)\cdot(1,0)=(0,0)$ випливає, що
\begin{align*}
  (i,i)&=(i,j)\cdot(j,i)=\\
       & =(0,1)\varphi\cdot(1,0)\varphi=\\
       & =((0,1)\cdot(1,0))\varphi=\\
       & =(0,0)\varphi=\\
       & =(0,0),
\end{align*}
а отже, $i=0$.

\smallskip

Приймемо $(0,1)\varphi=(0,k)$. Оскільки $(0,1)$ i $(1,0)$ --- твірні елементи біциклічного моноїда $\boldsymbol{B}_{\omega}$, то для довільного елемента $(m,n)\in\boldsymbol{B}_{\omega}$ маємо, що
\begin{align*}
  (m,n)\varphi &=((m,0)\cdot(0,n))\varphi=\\
               &=((1,0)^m\cdot(0,1)^n)\varphi= \\
               &=((1,0)^m)\varphi\cdot((0,1)^n)\varphi=\\
               &=((1,0)\varphi)^m\cdot((0,1)\varphi)^n=\\
               &=(k,0)^m\cdot(0,k)^n=\\
               &=(km,0)\cdot(0,kn)=\\
               &=(km,kn),
\end{align*}
звідки і випливає твердження леми.
\end{proof}

З леми~\ref{lemma-2.3} випливає, що
\begin{equation*}
  \varepsilon_{k_1[0]}\cdot\varepsilon_{k_2[0]}=\varepsilon_{k_1k_2[0]}=\varepsilon_{k_2k_1[0]}=\varepsilon_{k_2[0]}\cdot\varepsilon_{k_1[0]}
\end{equation*}
для довільних ендоморфізмів $\varepsilon_{k_1[0]}$ i $\varepsilon_{k_2[0]}$, $k_1,k_2\in\omega$, біциклічного моноїда $\boldsymbol{B}_{\omega}$, а отже, виконується така теорема.

\begin{theorem}\label{theorem-2.4}
Напівгрупа $\mathrm{\mathbf{End}}^1(\boldsymbol{B}_{\omega})$ ендоморфізмів біциклічного моноїда $\boldsymbol{B}_{\omega}$ ізо\-морф\-на напівгрупі $(\omega,*)$.
\end{theorem}

Рівності
\begin{equation}\label{eq-2.2}
  (k+a,a)^n=(kn+a,a) \qquad \hbox{i} \qquad (a,k+a)^n=(a,kn+a), \qquad a,k\in\omega, \quad n\in\mathbb{N},
\end{equation}
для елементів $(k+a,a), (a,k+a)$ біциклічного моноїда $\boldsymbol{B}_{\omega}$ доводяться методом математичної індукції.

\smallskip

Твердження~\ref{proposition-2.5} описує ендоморфізми біциклічного моноїда як напівгрупи.

\begin{proposition}\label{proposition-2.5}
Для довільного ендоморфізму $\varphi\colon \boldsymbol{B}_{\omega}\to \boldsymbol{B}_{\omega}$ біциклічного моноїда як напівгрупи існують такі  $k,a\in\omega$, що $\varphi=\varepsilon_{k[a]}$.
\end{proposition}

\begin{proof}
Оскільки $(0,0)$ --- ідемпотент біциклічного моноїда $\boldsymbol{B}_{\omega}$, то образ $(0,0)\varphi$ є ідемпотентом у $\boldsymbol{B}_{\omega}$.

\smallskip

У випадку $(0,0)\varphi=(0,0)$ твердження випливає з леми~\ref{lemma-2.3}. Тому надалі вважатимемо, що
\begin{equation*}
(0,0)\varphi=(a,a)\neq(0,0).
\end{equation*}
Оскільки $(0,1)\cdot(1,0)=(0,0)$, $(1,0)\cdot(0,1)=(1,1)$ i $(1,1)\preccurlyeq (0,0)$ в $\boldsymbol{B}_{\omega}$, то, врахувавши, що кожен гомоморфізм інверсних напівгруп зберігає природний частковий порядок, отримуємо, що $(1,1)\varphi\preccurlyeq(0,0)\varphi=(a,a)$ в $\boldsymbol{B}_{\omega}$. Якщо $(1,1)\varphi=(a,a)$, то за наслідком~1.32~\cite{Clifford-Preston-1961} ендоморфізм $\varphi$ є груповим, і оскільки всі підгрупи в біциклічному моноїді є тривіальними, то отримуємо, що $\varphi\colon \boldsymbol{B}_{\omega}\to \boldsymbol{B}_{\omega}$ --- анулюючий ендоморфізм, тобто $\varphi=\varepsilon_{0[a]}$ для деякого натурального числа $a$.

\smallskip

Надалі вважатимемо, що $(1,1)\varphi\neq(a,a)$. З означення природного часткового порядку на біциклічному моноїді $\boldsymbol{B}_{\omega}$ випливає, що існує таке натуральне число $k$, що
\begin{equation*}
(1,1)\varphi=(k+a,k+a).
\end{equation*}
Припустимо, що для твірних елементів $(0,1)$ i $(1,0)$ біциклічного моноїда $\boldsymbol{B}_{\omega}$ маємо, що $(0,1)\varphi=(i,j)$ i $(1,0)\varphi=(m,n)$ для деяких $i,j,m,n\in\omega$. Оскільки $\boldsymbol{B}_{\omega}$ --- інверсна напівгрупа та $(0,1)$ i $(1,0)$ --- інверсні елементи в $\boldsymbol{B}_{\omega}$, то
\begin{align*}
  (m,n)&=(1,0)\varphi=\\
   &=((0,1)^{-1})\varphi=\\
   &=((0,1)\varphi)^{-1}=\\
   &=(i,j)^{-1}=\\
   &=(j,i),
\end{align*}
а отже,
\begin{align*}
  (a,a) &=(0,0)\varphi=\\
   &=((0,1)\cdot(1,0))\varphi=\\
   &=(0,1)\varphi\cdot(1,0)\varphi=\\
   &=(i,j)\cdot(j,i)=\\
   &=(i,i),
\end{align*}
звідки випливає, що $i=a$. Також з аналогічних міркувань випливає, що
\begin{align*}
  (k+a,k+a) &=(1,1)\varphi=\\
  &=((1,0)\cdot(0,1))\varphi=\\
  &=(1,0)\varphi\cdot(0,1)\varphi=\\
  &=(j,a)\cdot(a,j)=\\
  &=(j,j),
\end{align*}
тобто $(0,1)\varphi=(a,k+a)$ i $(1,0)\varphi=(k+a,a)$. Тоді для довільного елемента $(m,n)$ біциклічного моноїда $\boldsymbol{B}_{\omega}$, використавши формули \eqref{eq-2.2}, отримуємо
\begin{align*}
  (m,n)\varphi&=((m,0)\cdot(0,n))\varphi=\\
   &=((1,0)^m\cdot(0,1)^n)\varphi= \\
   &=((1,0)^m)\varphi\cdot((0,1)^n)\varphi=\\
   &=((1,0)\varphi)^m\cdot((0,1)\varphi)^n=\\
   &=(k+a,a)^m\cdot(a,k+a)^n=\\
   &=(km+a,a)\cdot(a,kn+a)=\\
   &=(km+a,kn+a).
\end{align*}
Звідки випливає, що $\varphi=\varepsilon_{k[a]}$.
\end{proof}

Нехай $S$ i $T$ --- напівгрупи та $\mathfrak{f}\colon T\to \mathrm{\mathbf{End}}(S), \: t\mapsto f_t$ --- гомоморфізм. На декартовому добутку $S\times T$ означимо напівгрупову операцію так:
\begin{equation*}
  (s_1,t_1)\cdot (s_2,t_2)=(s_1\cdot (s_2)f_{t_1},t_1\cdot t_2).
\end{equation*}
Множина $S\times T$ з так визначеною напівгруповою операцією називається \emph{напівпрямим добутком} напівгрупи $S$ напівгрупою $T$ стосовно гомоморфізму $\mathfrak{f}$, та познача\-єть\-ся $S\rtimes_{\mathfrak{f}}T$ \cite{Lawson-1998}.

\smallskip

Означимо відображення $\mathfrak{f}\colon (\omega,*)\to \mathrm{\mathbf{End}}(\omega,+)$ за формулою $\mathfrak{f}(k)(n)=kn$. Очевидно, що ві\-доб\-раження $\mathfrak{f}$ є гомоморфізмом з напівгрупи $(\omega,*)$ у $\mathrm{\mathbf{End}}(\omega,+)$. Оскільки напівгрупа $(\omega,*)$ діє на напівгрупі $(\omega,+)$ ендоморфізмами, то на декартовому добутку $(\omega,+)\times(\omega,*)$ визначена напівгрупова операція
\begin{equation*}
  (a_1,k_1)\cdot(a_2,k_2)=(a_1+k_1a_2,k_1k_2)
\end{equation*}
напівпрямого добутку $(\omega,+)\rtimes_{\mathfrak{f}}(\omega,*)$ стосовно гомоморфізму $\mathfrak{f}$. Означимо ві\-доб\-ра\-жен\-ня  $\mathfrak{I}\colon \mathrm{\mathbf{End}}(\boldsymbol{B}_{\omega})\to(\omega,+)\rtimes_{\mathfrak{f}}(\omega,*)$ за формулою $\varepsilon_{k[a]}\mapsto (a,k)$. З леми \ref{lemma-2.1} і тверд\-жен\-ня \ref{proposition-2.5} випливає, що $\mathfrak{I}$ --- бієктивне відображення.

\smallskip

Отже, ми довели таку теорему.

\begin{theorem}\label{theorem-2.6}
Напівгрупа $\mathrm{\mathbf{End}}(\boldsymbol{B}_{\omega})$ ендоморфізмів біциклічної напівгрупи $\boldsymbol{B}_{\omega}$ ізоморфна напівпрямому добутку $(\omega,+)\rtimes_{\mathfrak{f}}(\omega,*)$.
\end{theorem}

\section{\textbf{Напівгрупа ендоморфізмів розширеної біциклічної напівгрупи}}\label{section-3}

Аналогічно, як і у випадку біциклічного моноїда $\boldsymbol{B}_{\omega}$ для довільних $k\in\omega$ i $a\in\mathbb{Z}$ означимо відображення $\varepsilon_{k[a]}\colon \boldsymbol{B}_{\mathbb{Z}}\to \boldsymbol{B}_{\mathbb{Z}}$ за формулою
\begin{equation}\label{eq-3.1}
  (m,n)\varepsilon_{k[a]}=(km+a,kn+a), \qquad m,n\in \mathbb{Z}.
\end{equation}

\smallskip

Доведення леми~\ref{lemma-3.1} аналогічне лемі~\ref{lemma-2.1}.

\begin{lemma}\label{lemma-3.1}
Для довільних $k\in\omega$ i $a\in\mathbb{Z}$ відображення $\varepsilon_{k[a]}\colon \boldsymbol{B}_{\mathbb{Z}}\to \boldsymbol{B}_{\mathbb{Z}}$, означене за формулою \eqref{eq-3.1}, є ендоморфізмом узагальненої біциклічної напівгрупи.
\end{lemma}

Рівності
\begin{equation}\label{eq-3.2}
  (k+a,a)^n=(kn+a,a) \qquad \hbox{i} \qquad (a,k+a)^n=(a,kn+a), \qquad a\in \mathbb{Z}, k,n\in\mathbb{N},
\end{equation}
для елементів $(k+a,a)$ i $(a,k+a)$ розширеної біциклічної напівгрупи $\boldsymbol{B}_{\mathbb{Z}}$, аналогічно як і формули \eqref{eq-2.2},  доводяться методом математичної індукції.

Піднапівгрупи розширеної біциклічної напівгрупи, які ізоморфні біциклічній напівгрупі, описано в лемі~\ref{lemma-3.2}.

\begin{lemma}\label{lemma-3.2}
Нехай $S$~--- піднапівгрупа розширеної біциклічної напівгрупи, яка ізоморфна біциклічній напівгрупі $\boldsymbol{B}_{\omega}$. Тоді
\begin{equation*}
S=\{(km+a,kn+a)\colon m,n\in\omega\}
\end{equation*}
для деяких $a\in\mathbb{Z}$ i $k\in\mathbb{N}$, причому ізоморфне занурення $\mathfrak{I}\colon \boldsymbol{B}_{\omega}\to\boldsymbol{B}_{\mathbb{Z}}$ визначається за формулою $(m,n)\mapsto(km+a,kn+a)$.
\end{lemma}

\begin{proof}
Нехай $\varphi\colon \boldsymbol{B}_{\omega}\to \boldsymbol{B}_{\mathbb{Z}}$~--- ізоморфне занурення біциклічної напівгрупи $\boldsymbol{B}_{\omega}$ у розширену біциклічну напівгрупу $\boldsymbol{B}_{\mathbb{Z}}$. Тоді $(0,0)\varphi=(a,a)$ для деякого цілого числа $a$. Оскільки $(1,1)\preccurlyeq(0,0)$ у $\boldsymbol{B}_{\omega}$, то $(1,1)\varphi\preccurlyeq(0,0)\varphi$ у $\boldsymbol{B}_{\mathbb{Z}}$, а отже за тверд\-женням~2.1$(i)$ з \cite{Fihel-Gutik-2011} існує натуральне таке число $k$, що $(1,1)\varphi=(k+a,k+a)$. Ми стверджуємо, що $(0,1)\varphi=(a,k+a)$ i $(1,0)\varphi=(k+a,a)$. Справді, припустимо, що для твірних елементів $(0,1)$ i $(1,0)$ біциклічного моноїда $\boldsymbol{B}_{\omega}$ маємо, що $(0,1)\varphi=(i,j)$ i $(1,0)\varphi=(m,n)$ для деяких $i,j,m,n\in\mathbb{Z}$. Оскільки $\boldsymbol{B}_{\omega}$ --- інверсна напівгрупа та $(0,1)$ i $(1,0)$ --- інверсні елементи в $\boldsymbol{B}_{\omega}$, то з того, що $\varphi\colon \boldsymbol{B}_{\omega}\to \boldsymbol{B}_{\mathbb{Z}}$~--- гомоморфізм інверсних напівгруп, випливає, що
\begin{align*}
  (m,n)& =(1,0)\varphi=\\
   & =((0,1)^{-1})\varphi=\\
   & =((0,1)\varphi)^{-1}=\\
   & =(i,j)^{-1}=(j,i),
\end{align*}
а отже,
\begin{align*}
  (a,a) &=(0,0)\varphi= \\
   & =((0,1)\cdot(1,0))\varphi= \\
   & =(0,1)\varphi\cdot(1,0)\varphi=\\
   & =(i,j)\cdot(j,i)=\\
   & =(i,i),
\end{align*}
звідки випливає, що $i=a$. Також з аналогічних міркувань випливає, що
\begin{align*}
  (k+a,k+a)& =(1,1)\varphi=\\
    & =((1,0)\cdot(0,1))\varphi=\\
    & =(1,0)\varphi\cdot(0,1)\varphi=\\
    & =(j,a)\cdot(a,j)=(j,j),
\end{align*}
тобто $(0,1)\varphi=(a,k+a)$ i $(1,0)\varphi=(k+a,a)$. Тоді для довільного елемента $(m,n)$ розширеної біциклічної напівгрупи  $\boldsymbol{B}_{\mathbb{Z}}$, використавши формули \eqref{eq-3.2}, отримуємо
\begin{align*}
  (m,n)\varphi&=((m,0)\cdot(0,n))\varphi=\\
   &=((1,0)^m\cdot(0,1)^n)\varphi= \\
   &=((1,0)^m)\varphi\cdot((0,1)^n)\varphi=\\
   &=((1,0)\varphi)^m\cdot((0,1)\varphi)^n=\\
   &=(k+a,a)^m\cdot(a,k+a)^n=\\
   &=(km+a,a)\cdot(a,kn+a)=\\
   &=(km+a,kn+a),
\end{align*}
що і завершує доведення леми.
\end{proof}

Ендоморфізм $\varphi\colon \boldsymbol{B}_{\mathbb{Z}}\to \boldsymbol{B}_{\mathbb{Z}}$ називається \emph{$(0,0)$-ендоморфізмом}, якщо
\begin{equation*}
(0,0)\varphi=(0,0).
\end{equation*}

\begin{lemma}\label{lemma-3.3}
Для довільного $(0,0)$-ендоморфізму $\varphi\colon \boldsymbol{B}_{\mathbb{Z}}\to \boldsymbol{B}_{\mathbb{Z}}$ розширеної бі\-цик\-лічної напівгрупи існує таке число $k\in\omega$, що $\varphi=\varepsilon_{k[0]}$.
\end{lemma}

\begin{proof}
Припустимо, що $(1,1)\varphi=(k,k)$.

\smallskip

Якщо $k=0$, то $(1,1)\varphi=(0,0)$. За твердженням~2.2 з \cite{Fihel-Gutik-2011} кожний неін'єктивний гомо\-морф\-ний образ узагальненої біциклічної напівгрупи $\boldsymbol{B}_{\mathbb{Z}}$ є циклічною групою, і оскільки за тверд\-жен\-ням~2.1$(iv)$ \cite{Fihel-Gutik-2011} усі максимальні підгрупи в $\boldsymbol{B}_{\mathbb{Z}}$ є тривіальними, то ендоморфізм $\varphi\colon \boldsymbol{B}_{\mathbb{Z}}\to \boldsymbol{B}_{\mathbb{Z}}$ анулюючий.

\smallskip

Припустимо, що $k\neq 0$. Оскільки $(1,1)\preccurlyeq(0,0)$ у напівгрупі $\boldsymbol{B}_{\mathbb{Z}}$, то $(k,k)\preccurlyeq(0,0)$, і тоді з  тверд\-жен\-ня~2.1$(i)$ \cite{Fihel-Gutik-2011} випливає, що $k>0$. Врахувавши, що множина
\begin{equation*}
\boldsymbol{B}_{\mathbb{Z}}[0]=\{(m,n)\colon m,n\in\omega\}
\end{equation*}
з індукованою напівгруповою операцією з $\boldsymbol{B}_{\mathbb{Z}}$ ізоморфна біциклічній напівгрупі, то з леми~\ref{lemma-2.3} випливає, що $(m,n)\varphi=(km,kn)$ для всіх $m,n\in\omega$.

\smallskip

За лемою~\ref{lemma-3.2} для довільного цілого числа $t$ множина
\begin{equation*}
\boldsymbol{B}_{\mathbb{Z}}[t]=\{(m,n)\colon m,n\geqslant t\}\subseteq \boldsymbol{B}_{\mathbb{Z}}
\end{equation*}
з індукованою з $\boldsymbol{B}_{\mathbb{Z}}$ напівгруповою операцією ізоморфна біциклічній напівгрупі стосовно відображення $(m+t,n+t)\mapsto(m,n)$.

\smallskip

Далі доведемо твердження леми методом математичної індукції. З наведеного вище випливає, що у випадку $p=0$ маємо, що $(m,n)\varphi=(km,kn)$ для всіх цілих $m,n\geqslant 0$. Припустимо, що з того, що виконується рівність $(m,n)\varphi=(km,kn)$ для всіх цілих $m,n\geqslant -p$, де $p$~--- деяке натуральне число, випливає, що ця рівність виконується для всіх цілих $m,n\geqslant -(p+1)$. Нехай
\begin{equation*}
(-(p+1),-(p+1))\varphi=(s,s)
\end{equation*}
для деякого цілого числа $s$. Оскільки
\begin{equation*}
(-p,-p)\preccurlyeq(-(p+1),-(p+1))
\end{equation*}
в $\boldsymbol{B}_{\mathbb{Z}}$ і кожен гомоморфізм інверсних напівгруп зберігає природний частковий порядок, то $(-kp,-kp)\preccurlyeq (s,s)$, і врахувавши, що $(-kp,-kp)\neq(s,s)$, то з тверд\-жен\-ня~2.1$(i)$ \cite{Fihel-Gutik-2011} випливає, що $s<-kp$. Тоді з рівностей
\begin{equation*}
(-p-1,-p)\cdot(-p,-p-1)=(-p-1,-p-1)
\end{equation*}
i
\begin{equation*}
(-p,-p-1)\cdot(-p-1,-p)=(-p,-p)
\end{equation*}
випливає, що
\begin{equation*}
  (-p-1,-p)\varphi=(s,-kp)  \qquad \hbox{i} \qquad (-p,-p-1)\varphi=(-kp, s).
\end{equation*}
Оскільки
\begin{equation*}
(-p,-p)\cdot (-p-1,-p)=(-p,-p+1),
\end{equation*}
то з припущення індукції та нерівності $s<-kp$ випливає, що
\begin{align*}
  (-kp,-kp+k)&=(-p,-p+1)\varphi= \\
   &=((-p,-p)\cdot(-p-1,-p))\varphi= \\
   &=(-p,-p)\varphi\cdot(-p-1,-p)\varphi= \\
   &=(-kp,-kp)\cdot(s,-kp)= \\
   &=(-kp+s-s,-kp-kp-s)= \\
   &=(-kp,-kp-kp-s),
\end{align*}
а отже, $-kp+k=-kp-kp-s$. Звідки випливає, що $s=-k(p+1)$. Тоді для довільного цілого числа $q<p+1$ маємо, що
\begin{equation*}
  (-p-1,-q)=(-p-1,-p)\cdot(-p,-q)
\end{equation*}
i
\begin{equation*}
  (-q,-p-1)=(-q,-p)\cdot (-p,-p-1),
\end{equation*}
а отже,
\begin{equation*}
\begin{split}
  (-p-1,-q)\varphi& =(-p-1,-p)\varphi\cdot(-p,-q)\varphi=\\
    & =(-kp-k,-kp)\cdot(-kp,-kq)=\\
    &=(-kp-k,-kq)
\end{split}
\end{equation*}
i
\begin{equation*}
\begin{split}
  (-q,-p-1)\varphi& =(-q,-p)\varphi\cdot (-p,-p-1)\varphi=\\
    & =(-kq,-kp)\cdot (-kp,-kp-k)=\\
    &=(-kq,-kp-k).
\end{split}
\end{equation*}
Звідси випливає, що $(m,n)\varphi=(km,kn)$ для всіх цілих $m,n\geqslant -p-1$, що і завершує доведення леми.
\end{proof}

З леми~\ref{lemma-3.3} випливає, що
\begin{equation*}
  \varepsilon_{k_1[0]}\cdot\varepsilon_{k_2[0]}=\varepsilon_{k_1k_2[0]}=\varepsilon_{k_2k_1[0]}=\varepsilon_{k_2[0]}\cdot\varepsilon_{k_1[0]}
\end{equation*}
для довільних $(0,0)$-ендоморфізмів $\varepsilon_{k_1[0]}$ i $\varepsilon_{k_2[0]}$, $k_1,k_2\in\omega$, розширеної біцик\-лічної на\-пів\-групи $\boldsymbol{B}_{\mathbb{Z}}$, а отже, виконується така теорема.

\begin{theorem}\label{theorem-3.4}
Напівгрупа $(0,0)$-ендоморфізмів розширеної біциклічної напівгрупи $\boldsymbol{B}_{\mathbb{Z}}$  ізоморфна напівгрупі $(\omega,*)$.
\end{theorem}

Твердження~\ref{proposition-3.5} описує ендоморфізми розширеної біциклічної напівгрупи.

\begin{proposition}\label{proposition-3.5}
Для довільного ендоморфізму $\varphi\colon \boldsymbol{B}_{\mathbb{Z}}\to \boldsymbol{B}_{\mathbb{Z}}$ розширеної біциклічної напівгрупи існують такі  $k\in\omega$ і $a\in\mathbb{Z}$, що $\varphi=\varepsilon_{k[a]}$.
\end{proposition}

\begin{proof}
Нехай $\varphi$~--- довільний ендоморфізм розширеної біциклічної напівгрупи  $\boldsymbol{B}_{\mathbb{Z}}$. Тоді  $(0,0)\varphi=(a,a)$ для деякого $a\in\mathbb{Z}$. За теоремою 1 \cite{Gutik-Maksymyk-2017} ендоморфізм $\varepsilon_{1[-a]}$ є елементом групи одиниць напівгрупи ендоморфізмів розширеної біциклічної напівгрупи $\boldsymbol{B}_{\mathbb{Z}}$. З визначення відображення $\varepsilon_{1[-a]}$ випливає, що
\begin{equation*}
(0,0)\varphi\varepsilon_{1[-a]}=(a,a)\varepsilon_{1[-a]}=(0,0),
\end{equation*}
а отже, відображення $\varphi\varepsilon_{1[-a]}$ є $(0,0)$-ендоморфізмом розширеної біциклічної напівгрупи. За лемою \ref{lemma-3.3} існує таке число $k\in\omega$, що $\varphi\varepsilon_{1[-a]}=\varepsilon_{k[0]}$. Оскільки $\varepsilon_{1[-a]}\varepsilon_{1[a]}=\varepsilon_{1[0]}$ --- тотожний автоморфізм розширеної біциклічної напівгрупи  $\boldsymbol{B}_{\mathbb{Z}}$, то
\begin{equation*}
\varphi=\varphi\varepsilon_{1[-a]}\varepsilon_{1[a]}=\varepsilon_{k[0]}\varepsilon_{1[a]}=\varepsilon_{k[a]},
\end{equation*}
що і завершує доведення твердження.
\end{proof}

Аналогічно, як і у випадку біциклічної напівгрупи, означимо відображення $\mathfrak{f}\colon (\omega,*)\to \mathrm{\mathbf{End}}(\mathbb{Z}(+))$ за формулою $\mathfrak{f}(k)(n)=kn$. Очевидно, що ві\-доб\-раження $\mathfrak{f}$ є гомоморфізмом з напівгрупи $(\omega,*)$ у напівгрупу $\mathrm{\mathbf{End}}(\mathbb{Z}(+))$. Оскільки напівгрупа $(\omega,*)$ діє на групі $\mathbb{Z}(+)$ ендоморфізмами, то на декартовому добутку $\mathbb{Z}(+)\times(\omega,*)$ визначена напівгрупова операція
\begin{equation*}
  (a_1,k_1)\cdot(a_2,k_2)=(a_1+k_1a_2,k_1k_2)
\end{equation*}
напівпрямого добутку $\mathbb{Z}(+)\rtimes_{\mathfrak{f}}(\omega,*)$ стосовно гомоморфізму $\mathfrak{f}$. Означимо ві\-доб\-ра\-жен\-ня  $\mathfrak{I}\colon \mathrm{\mathbf{End}}(\boldsymbol{B}_{\mathbb{Z}})\to\mathbb{Z}(+)\rtimes_{\mathfrak{f}}(\omega,*)$ за формулою $\varepsilon_{k[a]}\mapsto (a,k)$. З леми \ref{lemma-3.1} і твердження~\ref{proposition-3.5} випливає, що $\mathfrak{I}$ --- бієктивне відображення. Отже, ми довели таку теорему.

\begin{theorem}\label{theorem-3.6}
Напівгрупа $\mathrm{\mathbf{End}}(\boldsymbol{B}_{\mathbb{Z}})$ ендоморфізмів розширеної біциклічної напівгрупи  $\boldsymbol{B}_{\mathbb{Z}}$ ізоморфна напівпрямому добутку $\mathbb{Z}(+)\rtimes_{\mathfrak{f}}(\omega,*)$.
\end{theorem}

\smallskip

\section*{\textbf{Подяка}}

Автори висловлюють щиру подяку Т. Банаху, О. Равському та  рецензентові за цінні поради та зауваження.


\end{document}